\documentclass[11pt,twoside]{article}

\usepackage{latexsym}
\usepackage{amsmath}
\usepackage{amsthm}
\usepackage{amssymb}
\usepackage{vmargin}
\usepackage{amscd}
\usepackage{stmaryrd}
\usepackage{euscript}
\usepackage{mathrsfs}
\usepackage{amscd}
\usepackage[all]{xy}
\usepackage{xr}
\DeclareMathAlphabet{\mathpzc}{OT1}{pzc}{m}{it}

\setmargins{32mm}{20mm}{14.6cm}{22cm}{1cm}{1cm}{1cm}{1cm}

\def\s|{\,|\,}
\def\ten{\otimes}

\def\CC{\mathrm{C}}

\def\SS{\mathrm{S}}

\def\Z{\mathbb{Z}}
\def\Q{\mathbb{Q}}
\def\R{\mathbb{R}}
\def\Cx{\mathbb{C}}

\def\vv{\mathbb{V}}
\def\ww{\mathbb{W}}
\def\Bu{\mathbb{B}}

\def\bB{\mathbb{B}}

\def\bD{\mathbb{D}}

\def\bF{\mathbb{F}}

\def\bH{\mathbb{H}}

\def\C{\mathcal{C}}

\def\cE{\mathcal{E}}

\def\cH{\mathcal{H}}

\def\cL{\mathcal{L}}

\def\cQ{\mathcal{Q}}

\def\cT{\mathcal{T}}

\def\cW{\mathcal{W}}

\def\F{\mathscr{F}}

\def\sA{\mathscr{A}}

\def\sC{\mathscr{C}}

\def\sF{\mathscr{F}}

\def\sP{\mathscr{P}}

\def\sZ{\mathscr{Z}}

\def\r{\mathfrak{R}}

\def\H{\mathrm{H}}
\def\z{\mathrm{Z}}
\def\b{\mathrm{B}}

\def\m{\mathfrak{m}}
\def\n{\mathfrak{n}}
\def\g{\mathfrak{g}}

\def\bu{\ad \Bu}

\def\ENd{\mathscr{E}\!\mathit{nd}}

\def\aut{\mathscr{A}\!\mathit{ut}}

\def\Hom{\mathrm{Hom}}

\def\Aut{\mathrm{Aut}}

\def\im{\mathrm{Im\,}}

\def\Set{\mathrm{Set}}

\def\ad{\mathrm{ad}}
\def\d{\tilde{d}}
\def\<{\langle}
\def\>{\rangle}
\def\Lim{\varprojlim}
\def\into{\hookrightarrow}
\def\onto{\twoheadrightarrow}
\def\xra{\xrightarrow}

\def\by{\times}
\def\mcl{\mathrm{MC}_L}
\def\mc{\mathrm{MC}}
\def\gl{\mathrm{G}_L}

\def\defl{\mathrm{Def}_L}

\def\defn{\mathrm{Def}_N}
\def\ddef{\mathrm{Def}}

\def\Einf{\mathrm{E}_{\mathrm{C}^{\infty}}}

\def\GL{\mathrm{GL}}

\def\etale{\acute{\mathrm{e}}\mathrm{tale}}

\def\diag{\mathrm{diag}}

\def\sigmah{\sigma_{\mathrm{h}}}
\def\sigmav{\sigma_{\mathrm{v}}}

\def\pd{\partial}
\def\dc{d^{\mathrm{c}}}
\def\half{\frac{1}{2}}

\newtheorem{theorem}{Theorem}[section]

\newtheorem{corollary}[theorem]{Corollary}
\newtheorem{lemma}[theorem]{Lemma}
\newtheorem{theorem*}{Theorem}
\newtheorem{proposition*}[theorem*]{Proposition}
\newtheorem{corollary*}[theorem*]{Corollary}
\newtheorem{lemma*}[theorem*]{Lemma}

\theoremstyle{definition}
\newtheorem{definition}[theorem]{Definition}

\newtheorem{remark}[theorem]{Remark}
\newtheorem{remarks}[theorem]{Remarks}

\newtheorem{definition*}[theorem*]{Definition}
\newtheorem{example*}[theorem*]{Example}
\newtheorem{examples*}[theorem*]{Examples}
\newtheorem{remark*}[theorem*]{Remark}
\newtheorem{remarks*}[theorem*]{Remarks}

\begin{document}
\title{The deformation theory of representations of the fundamental group of a smooth variety}
\author{J.P.Pridham\thanks{The author is supported by Trinity College, Cambridge and the Isle of Man Department of Education.}}
\maketitle

%
\tableofcontents

\section*{Introduction}
\addcontentsline{toc}{section}{Introduction}

In \cite{sdc}, I developed the theory of Simplicial Deformation Complexes (SDCs), proposed as an alternative to differential graded Lie algebras (DGLAs). They can be constructed for a whole range of deformation problems and capture more information than just the deformation functor. On a more practical level, DGLAs can be extremely useful; they enable Goldman and Millson in \cite{GM} to show, with certain restrictions, that the hull of the functor describing deformations of a real representation of the fundamental group of a compact K\"ahler manifold is defined by homogeneous quadratic equations.

Even for the problem considered in \cite{GM}, where  a perfectly satisfactory governing DGLA was constructed, SDCs can be helpful. Goldman and Millson proceed via a chain of groupoid equivalences from the deformation groupoid to the groupoid associated to the DGLA. However, it is possible to write down an SDC immediately, passing via a natural chain of quasi-isomorphisms to the DGLA in \cite{GM}. This substantially shortens the reasoning, since it is no longer necessary to find topological interpretations of each intermediate object and functor, and quasi-isomorphisms are quicker to establish than groupoid equivalences.
 
\cite{GM} was motivated by \cite{DGMS}, which shows that the homotopy type of a compact K\"ahler manifold is a formal consequence of its cohomology, the idea being to replace statements about DGAs with those about DGLAs.  
In \cite{Morgan}, Morgan  proved an analogue of the results in \cite{DGMS}  for smooth complex varieties. This indicates that the  fundamental group of such a variety should have similar properties to those established in \cite{GM} for compact K\"ahler manifolds.  

In Section \ref{toprep}, I use the theory of SDCs to give a shorter proof of Goldman and Millson's result. I also prove that the hull of the deformation functor of representations of the (topological) fundamental group of a smooth (non-proper) complex variety has a mixed Hodge structure. Using the weight restrictions established in Hodge II (\cite{Hodge2}), it follows that the hull is defined by equations of degree at most four. Neither of these results employs the theory of SDCs in an essential way. However, simplicial methods are vital for the final result of the section --- that there is a mixed Hodge structure on the hull, even when the variety is not smooth. This has few consequences, since the weight restrictions of Hodge III (\cite{Hodge3}) give no bound on the degree of the defining equations.

The results in \cite{DGMS} were inspired by the yoga of weights, motivated by the Weil Conjectures.  With this in mind, it is reasonable to expect the characteristic zero results to have finite characteristic analogues. Section \ref{algrep} is concerned with proving results on the structure of the hull of the functor of deformations of a continuous $l$-adic representation of the algebraic fundamental group of a smooth variety in finite characteristic. This is done by studying the behaviour of the Frobenius action on the hull. As the obstruction maps are Frobenius equivariant, the SDC allows us to describe the hull in terms of cohomology groups. The weight restrictions in Weil II (\cite{Weil2}) then imply that the hull is quadratic for smooth proper varieties, and defined by equations of degree at most four for smooth, non-proper varieties.

Although the two sections of this paper prove similar results by using the same underlying philosophy, namely that weights determine the structure of the hull, they are logically independent. 

\section{Representations of the Topological Fundamental Group}\label{toprep}

Fix a connected topological space  $X$ (sufficiently nice to have a universal covering space) and a  point $x \in X$. Denote $\pi_1(X,x)$ by $\Gamma$. Throughout, $A$ will denote a ring in $\C_{\R}$.
Let $G$ be a  real algebraic Lie group, and fix a representation $\rho_0: \Gamma \to G(\R)$.

\begin{definition} Given a  representation $\rho: \Gamma \to G(A)$, let
$$
\Bu_{\rho}:=\F_{\tilde{X}}\by_{\Gamma, \rho}G(A)=(\pi_*G(A))^{\Gamma,\rho},
$$
where $\F_{\tilde{X}}$ is the sheaf on $X$ whose espace \'etal\'e is $\tilde{X}$, and  $ \tilde{X}\xra{\pi} X$ is the universal covering space of $X$.
 
\end{definition}

\begin{lemma} The map
\begin{eqnarray*}
\r_A(\rho_0) &\to& \mathfrak{B}_A(\Bu_0)\\
\rho &\mapsto& \Bu_{\rho}
\end{eqnarray*}
is an equivalence of categories,
where $\r_A(\rho_0)$ has as  objects  representations 
$$
\rho \equiv \rho_0 \mod \m_A,
$$ 
while $\mathfrak{B}_A(\Bu_0)$ has as objects principal $G(A)$-sheaves  
$$
\Bu \equiv \Bu_0 \mod \m_A
$$ on $X$. In both cases morphisms must be the identity mod $ \m_A$. 
\begin{proof}
This is essentially the same as the opening pages of \cite{De}.
It suffices to construct a quasi-inverse. Given a principal $G(A)$-sheaf $\Bu$ on $X$, $\pi^{-1}\Bu$ will be a principal $G(A)$-sheaf on $\tilde{X}$, so will be a constant sheaf, since $\tilde{X}$ is simply connected. Fix a point $y \in \tilde{X}$ above $x$. We have, for each $\gamma \in \Gamma$,  isomorphisms
$$
\Bu_x \cong (\pi^{-1}\Bu)_y \cong (\pi^{-1}\Bu)_{\gamma y} \cong \Bu_x
$$
of principal $G(A)$-spaces, the second isomorphism arising from the constancy of the sheaf. The composite isomorphism can be taken to be right multiplication by $\rho(\gamma)$, thus defining a representation $\rho$.
\end{proof} 
\end{lemma}

Thus, as in \cite{sdc} Section \ref{sdc-smoothgroup}, this deformation problem is governed by the SDC
$$
E^n=\Hom(u^{-1}\Bu_0,(u^{-1}u_*)^n u^{-1}\Bu_0)^{G(\R)}_{u^{-1}\alpha^n}.
$$
The cohomology groups are isomorphic to 
$$
\H^*(X,\bu_0),
$$ 
where $\bu_0$ is the tangent space of $\aut(\Bu_0)$ at the identity. Equivalently,
$$
\bu_0=\g \by^{G(\R)}\Bu_0= (\pi_* \g)^{\Gamma,\rho_0},
$$
for $\g$ the tangent space of $G(\R)$ at the identity, and $G(\R)$ acting via the adjoint action,
so the cohomology groups are
$$
\H^*(X,(\pi_*\g)^{\Gamma,\rho_0}).
$$

Recall from \cite{sdc} Section \ref{sdc-sdcdgla} that there is a functor $\cE$ from DGLAs to SDCs. In particular, if $L=\g$, a Lie algebra concentrated in degree $0$, then 
$$
\cE^i(\g)=\ker(G(A) \to G(\R))=\exp(\g \ten \m_A),
$$ 
where $G$ is the corresponding Lie group, with all $\sigma, \pd$ the identity, and $*$ given by the Alexander-Whitney cup product
$$
e*f = (\pd^{m+n}\ldots \pd^{m+2}\pd^{m+1}e)\cdot (\pd^0)^m f,
$$
 for $e \in \cE^m,\,f \in \cE^n$.

Provided that $\H^1(X,\bu_0)$ is finite dimensional, the deformation functor will have a hull.

Observe the quasi-isomorphism of SDCs
\begin{eqnarray*}\label{bungp}
\Gamma(X', (u^{-1}u_*)^n u^{-1} \exp(\bu_0)) &\to& \Hom(u^{-1}(\Bu_0),(u^{-1}u_*)^n u^{-1}(\Bu_0))^{G(\R)}_{u^{-1}\alpha^n},\\
g  &\mapsto& (b \mapsto g\cdot b),
\end{eqnarray*} 
where  the SDC structure of the term on the left is given by the Alexander-Whitney cup product.

Now, write 
$$
\CC^n(X,\sF):= \Gamma(X', (u^{-1}u_*)^n u^{-1}\sF),
$$
 for sheaves $\sF$ on $X$, and
$$
\sC^n(\sF):=(u_* u^{-1})^{n+1}\sF,
$$
so that $\CC^n(X,\sF)=\Gamma(X, \sC^n(\sF))$.

Now assume that $X$ is a differentiable manifold, and let $\sA^n$ be the sheaf of real-valued $\mathrm{C}^{\infty}$ $n$-forms on $X$.

We have the following quasi-isomorphisms of SDCs:
$$
\begin{matrix}
\CC^n(X,\exp(\bu_0)) \\
\|	\\
\CC^n(X,\cE(\bu_0)^n) \\
\downarrow \\
\CC^n(X,\cE(\bu_0 \ten \sA^{\bullet})^n)\\
\uparrow \\
\cE(\Gamma(X,\bu_0 \ten \sA^{\bullet}))^n,
\end{matrix}
$$
since on cohomology we have:
$$
\begin{matrix}
\H^n(X,\bu_0) \\
\downarrow \\
\bH^n(X,\bu_0 \ten \sA^{\bullet})\\
\uparrow \\
\H^n(\Gamma(X,\bu_0 \ten \sA^{\bullet})),
\end{matrix}
$$
the first quasi-isomorphism arising because $\R \to \sA^{\bullet}$ is a resolution, and
the second  from the flabbiness of $\sA^{\bullet}$.

Hence $\r_A(\rho_0)$ is governed by the DGLA
$$
\Gamma(X,\bu_0 \ten \sA^{\bullet}).
$$

\subsection{Compact K\"ahlerian manifolds}\label{compactkahler}
This section is a reworking of \cite{GM}. Let $X$ be a  compact complex K\"ahlerian manifold.

If $\ad(\rho_0)(\Gamma) \subset K$, for some compact subgroup $K \le \GL(\g)$, then we can define $K$-invariant Haar measure on $\g$, and thus a global positive definite inner product on $\bu_0$, so it becomes an orthogonal local system.

 On $\sA^*$, we have not only the operator $d$, but also
$$
\dc:=J^{-1}dJ=\sqrt{-1}(\bar{\pd}-\pd),
$$
where $J$ is the integrable almost-complex structure on $X$. Define
$$
\sZ^{\bullet}_{\dc}:= \ker (\dc: \sA^{\bullet} \to \sA^{\bullet}).
$$

 We now apply Hodge theory for orthogonal local systems, in the form of the $d \dc$ Lemma (identical to \cite{DGMS} 5.11):
\begin{lemma} 
On $\Gamma(X, \sA^n \ten \bu_0)$, 
\begin{eqnarray*}
\ker(\dc)\cap \im(d)=\im(d \dc),\\
\ker(\d)\cap \im(\dc)=\im(d \dc).
\end{eqnarray*}
\end{lemma}

\begin{remark}
It is worth noting that the $d\dc$ lemma holds on any manifold which can be blown up to a compact K\"ahler manifold, hence on any compact Moishezon manifold, so we may take $X$ to be of this form. Of course, this gives no more information about fundamental groups, since blowing-up does not change the fundamental group.
\end{remark}

We thus obtain quasi-isomorphisms of DGLAs:
$$
\begin{matrix}
\Gamma(X,\bu_0 \ten \sA^{\bullet})\\
\uparrow \\
\Gamma(X,\bu_0 \ten \sZ_{\dc}^{\bullet})\\
\downarrow \\
(\H_{\dc}^*(\bu_0\ten \sA^{\bullet}),d=0)\\
\| \\
(\H^*(X,\bu_0),0).
\end{matrix}
$$
To see that these are quasi-isomorphisms, we look at the cohomology groups:
$$
\begin{matrix}
\H^n(\Gamma(X,\bu_0 \ten \sA^{\bullet}))\\
\uparrow \\
\H^n(\Gamma(X,\bu_0 \ten \sZ_{\dc}^{\bullet},d))\\
\downarrow \\
(\H^n_{\dc}(\Gamma(X,\bu_0\ten \sA^{\bullet})).
\end{matrix}
$$

 The last two quasi-isomorphisms are a consequence of the $d\dc$ Lemma.

However,  the DGLA $(H^*(\bu_0),d=0)$  is formal, in the sense of \cite{Man}, so the hull of the deformation functor is given by homogeneous quadratic equations. Note that a hull $S$ must exist, since $\dim \H^1(\bu_0) < \infty$. Explicitly,
$$
S \cong \R[[\H^1(\bu_0)^{\vee}]]/(\check{\cup}\,\H^2(\bu_0)^{\vee}),
$$
where
$$
\check{\cup}\,:\H^2(\bu_0)^{\vee} \to \SS^2(\H^1(\bu_0)^{\vee})
$$
is dual to the cup product.

Note further that we have a Hodge decomposition on $S \ten \Cx$, where $S$ is the hull. Since $\m_S/\m_S^2$ is pure of weight $-1$, $\m_S^n/\m_S^{n+1}$ is pure of weight $-n$, so must have even dimension for $n$ odd. 

\subsection{Smooth complex varieties}\label{smoothcx}

In this case we have to use the mixed Hodge theory of \cite{Hodge2}, rather than pure Hodge theory.
Let $X$ be a smooth complex variety. Using \cite{Nagata} and \cite{Hironaka}, we may write $X=\bar{X} - D$, for $\bar{X}$ a complete smooth variety, and $D$ a divisor with normal crossings.
 Again, the governing DGLA is 
$$
\Gamma(X,\bu_0 \ten \sA^{\bullet}).
$$

If $\ad(\rho_0)(\Gamma) \subset K$, for some compact subgroup $K \le \GL(\g)$, then we can define $K$-invariant Haar measure on $\g$, and thus a global positive definite inner product on $\bu_0$, so that it becomes an orthogonal local system. 

Where \cite{GM} mimicked \cite{DGMS}, we must now mimic \cite{Morgan}, and will adopt the notation of that paper. 

Observe that \cite{Morgan} \S 2 works for twisted coefficients (in our case $\bu_0$), since $\bu_0$ is trivial on each simplex,  giving us a quasi-isomorphism of DGLAs
$$
\Einf(X, \bu_0) \to \Gamma(X,\bu_0 \ten \sA^{\bullet}),
$$
There is then a quasi-isomorphism
$$
\Einf(X, \bu_0)\ten \Cx \to (\Gamma(\bar{X},\widetilde{\bu_0}\ten \Omega_{\bar{X}}^{*}(\log D)\ten \sA_{\Cx}^{0,*}), \nabla),
$$
where $(\widetilde{\bu_0}, \nabla)$ is the holomorphic Deligne extension of $\bu_0$, as defined in \cite{De}.
The mixed Hodge theory of \cite{Timm} now proves that the latter quasi-isomorphism is a mixed Hodge diagram, in the sense of \cite{Morgan} 3.5 (except that we work with DGLAs rather than DGAs).

In \cite{Sullivan}, Sullivan defines the minimal model of a connected DGA. Replacing free graded algebras by free graded Lie algebras yields the same construction for a DGLA $L$ in non-negative degrees with $\H^0(L)=0$. The minimal model $M\to L$ will then in some sense be universal in the quasi-isomorphism class of $L$. If we only require quasi-isomorphisms in degrees $\ge 1$, then we can analogously construct a minimal model $M \to L$, for which automatically $\H^0(M)=0$.

Let $M$ denote the Sullivan minimal model of $ \Gamma(X,\bu_0 \ten \sA^{\bullet})$. By \cite{Morgan} Theorems 6.10 and 8.6, there is a mixed Hodge structure on $M$. By \cite{Morgan} Theorem 10.1, this induces a weight decomposition on $M$. 

\begin{theorem}\label{morthm}

The hull $S$ of the functor $\ddef_{\rho_0}$ is defined by equations of degree at most $4$. Moreover, there is a mixed Hodge structure on $S$, with $\m_S/\m_S^2$ of weights $-1$ and $-2$.

\begin{proof}
Choose a  decomposition 
$$
M^1=\H^1(M)\oplus C^1, \quad \z^2(M)= \cH^2 \oplus \b^1(M), \quad M^2= \z^2(M)\oplus C^2,
$$
respecting the filtrations $(F,W)$ on $M$ (this is straightforward, since $F,\bar{F},W$ give a bigraded decomposition of $M\ten \Cx$).

The hull of $\ddef_M$ is then given by the Kuranishi functor as defined in \cite{Man} \S 4, which will now be pro-represented by
$$
\R[[\H^1(M)^{\vee}]]/(f\cH^2(M)^{\vee}),
$$
where 
$$
f:\cH^2(M)^{\vee} \to (\H^1(M)^{\vee})^2 \triangleleft  \R[[\H^1(M)^{\vee}]]
$$ 
respects the filtrations. By \cite{Timm} \S 6, $\H^1(X,\bu_0)$ is of weights $1$ and $2$, while $\H^2(X,\bu_0)$ is of weights $2, 3$ and $4$. Thus $f$ corresponds to functions:
\begin{eqnarray*}
(W_{-1})^2 &\text{weight }-2 \\
(W_{-1})^3 +(W_{-1})(W_{-2}) &\text{weight }-3\\
(W_{-1})^4 +(W_{-1})^2(W_{-2}) +(W_{-2})^2 &\text{weight }-4,
\end{eqnarray*}
hence equations of degree at most $4$. 

Moreover, $f$ respects the Hodge structure, so we have a mixed Hodge structure on $S$. Note that the quotient map
$$
f:\H^2(M)^{\vee} \to (\H^1(M)^{\vee})^2/(\H^1(M)^{\vee})^3 \cong S^2(\H^1(M)^{\vee})
$$
is dual to half the cup product
\begin{eqnarray*}
\H^1(X,\bu_0) &\to&  \H^2(X,\bu_0),\\
a &\mapsto& \half a\cup a.
\end{eqnarray*}
\end{proof}
\end{theorem}

\begin{remark} \label{hodgedef}
It ought to be possible to prove this result without recourse to minimal models and the methods of \cite{Morgan}. The following approach should work:
\begin{enumerate}
\item Define a category $\mathrm{Hdg}(\C_{\R})$ of real Artinian algebras, with mixed Hodge structures on their maximal ideals. Show that Schlessinger's theorems carry over into this context (in particular, tangent and obstruction spaces will have mixed Hodge structures).

\item Given a filtration $F$ on a  complex  DGLA $L$, and an algebra $A \in \mathrm{Hdg}(\C_{\R})$, define 
\begin{eqnarray*}
F \mcl(A)&:=& \mcl(A\ten \Cx)\cap F^0(L^1\ten \m_A)\\
F \gl(A)&:=& \exp (d^{-1} F^0(L^1\ten \m_A))\\
F \defl(A)&:=& F \mcl(A)/F \gl(A).
\end{eqnarray*}
We have  similar definitions for a filtration $\bar{F}$ on a complex DGLA, and for a filtration $W$ on a real DGLA. 

\item Consider the DGLA 
$$
N=(\Gamma(\bar{X},\widetilde{\bu_0}\ten \Omega_{\bar{X}}^{*}(\log D)\ten \sA_{\Cx}^{0,*}), \nabla),
$$
defined as above. This has filtrations $W$ and $F$. We have quasi-isomorphisms (of $W$-filtered DGLAs)
$$
 N \to \Gamma(X,\bu_0 \ten \sA^{\bullet})\ten \Cx \leftarrow \bar{N}.
$$
Writing $L=\Gamma(X,\bu_0 \ten \sA^{\bullet})$, we can use these quasi-isomorphisms to identify the functors $\ddef_N, \ddef_{\bar{N}}$ and $\ddef_{L\ten \Cx}$.

Now, consider the functor on $\mathrm{Hdg}(\C_{\R})$ given by
$$
A \mapsto W \defl(A) \cap F\defn(A\ten \Cx) \cap \bar{F}\ddef_{\bar{N}}(A\ten \Cx),
$$
the intersection being as subsets of $\defl(A\ten \Cx)$. This functor should have a hull $S$, with tangent space $\H^1(X,\bu_0)$ and obstruction space $\H^2(X,\bu_0)$ (with the usual Hodge structures). It would remain only to show that $S$ (without its Hodge structures) is a hull for $\defl$ over $\C_{\R}$.
\end{enumerate}
\end{remark}

\subsection{Arbitrary complex varieties}\label{arbcx}

In this section, the methods of \cite{Hodge3} will be used. Let $X$ be a complex variety. As in \cite{Hodge3} 8.1.12, we take a simplicial resolution $Y_{\bullet} \to X$ of $X$ by smooth varieties, and compactify to obtain a smooth proper simplicial variety $\overline{Y}_{\!\!\bullet}$, and a divisor $D_{\bullet}$ of normal crossings in $\overline{Y}_{\!\!\bullet}$, such that $Y_{\bullet}=\overline{Y}_{\!\!\bullet} - D_{\bullet}$. 

Now, the SDC governing deformations of $\rho_0$ is
$$
\CC^n(X, \exp(\bu_0)).
$$

Cohomological descent implies that we have a quasi-isomorphism of SDCs:
$$
\CC^n(X, \exp(\bu_0)) \to \CC^n(Y_n, \exp(\bu_0)).
$$
We now proceed as in the previous section, obtaining quasi-isomorphisms

$$
\begin{matrix}
\CC^n(Y_n, \exp(\bu_0)) \\
\downarrow\\
 \CC^n(Y_n,\cE(\bu_0 \ten \sA_{Y_n}^{\bullet})^n) \\
\uparrow\\
 \cE(\Gamma(Y_n,\bu_0 \ten \sA_{Y_n}^{\bullet}))^n\\
\uparrow\\
\cE(\Einf(Y_n, \bu_0))^n .
\end{matrix}
$$

As in \cite{Morgan} \S 3, there is a quasi-isomorphism
$$
\Einf(Y_n, \bu_0)\ten \Cx \to (\Gamma(\overline{Y}_{\!\!n},\widetilde{\bu_0}\ten \Omega_{\overline{Y}_{\!\!n}}^{*}(\log D_n)\ten \sA_{\Cx}^{0,*}), \nabla),
$$
with $\nabla$ as in \cite{Timm}. We have a weight filtration $W$ on the former, and a compatible pair $(W,F)$ of filtrations on the latter.

In \cite{Hodge3} \S 7, Deligne defines filtrations $F,\delta(W,L)$ on the total complex of a bigraded cochain complex. Since we are working simplicially, we must, via the Eilenberg-Zilber Theorem, define the corresponding filtration on the diagonal complex of a bi-cosimplicial complex.

As in \cite{sdc} Section \ref{sdc-explicit}, we form the cosimplicial Lie algebras
$$
\mathpzc{e}(\Einf(Y_n, \bu_0))\ten \Cx \to \mathpzc{e}(\Gamma(\overline{Y}_{\!\!n},\widetilde{\bu_0}\ten \Omega_{\overline{Y}_{\!\!n}}^{*}(\log D_n)\ten \sA_{\Cx}^{0,*}), \nabla),
$$

Now, we have
$$
\mathpzc{e}(\Einf(Y_n, \bu_0))^n = \diag\mathpzc{e}(\Einf(Y_m, \bu_0))^n,
$$
and similarly for the other complex, on which we wish to define filtrations.
Given a cosimplicial complex $K$ of vector spaces, the Eilenberg-Zilber map
$$
(\diag K)^n= K^{nn} \xra{\nabla} \bigoplus_{p+q=n} K^{pq}  =(\mathrm{Tot} K)^n
$$
induces an isomorphism on cohomology. Moreover, if 
$$
\alpha \in N^n(\diag K)= \bigcap_{i=0,\ldots,n} \ker(\sigmah^i\sigmav^i: (\diag K)^n \to (\diag K)^{n-1}),
$$
then $\nabla(\alpha)_{pq} \in N^p_{\mathrm{h}}\cap N^q_{\mathrm{v}}$.

For a DGLA $L$, observe that $N(\mathpzc{e}(L))=L$. We therefore define
$$
\delta(W,L)_n(N(\diag \mathpzc{e}(\Einf(Y_m, \bu_0))^n)) = \nabla^{-1} \delta(W,L)_n (\mathrm{Tot} \Einf(Y_{\bullet}, \bu_0)^{\bullet}),
$$
the latter being defined as in \cite{Hodge3} \S 7. Similarly, we obtain filtrations $\delta(W,L)$ and $F$ on
$$
N(\mathpzc{e}(\Gamma(\overline{Y}_{\!\!n},\widetilde{\bu_0}\ten \Omega_{\overline{Y}_{\!\!n}}^{*}(\log D_n)\ten \sA_{\Cx}^{0,*}))).
$$

Now, given a cosimplicial complex $K^*$, and a filtration $\mathrm{Fil}$ on $N(K)$,
define
$$
\mathrm{Fil}_m(K^n)= \sum_{0 \le j_1 < \ldots < j_t< n}  (\sigma^{j_1}\ldots \sigma^{j_t})^{-1} \mathrm{Fil}_m(N^{n-t}(K)).
$$
The reasoning behind these constructions is that they preserve the filtrations induced on cohomology.

We have therefore defined filtrations $\delta(W,L),F$ (as applicable) on
$$
\mathpzc{e}(\Einf(Y_n, \bu_0))\ten \Cx \to \mathpzc{e}(\Gamma(\overline{Y}_{\!\!n},\widetilde{\bu_0}\ten \Omega_{\overline{Y}_{\!\!n}}^{*}(\log D_n)\ten \sA_{\Cx}^{0,*}), \nabla),
$$

Exponentiation induces corresponding filtrations on
$$
\cE(\Einf(Y_n, \bu_0)\ten \Cx) \to \cE(\Gamma(\overline{Y}_{\!\!n},\widetilde{\bu_0}\ten \Omega_{\overline{Y}_{\!\!n}}^{*}(\log D_n)\ten \sA_{\Cx}^{0,*}), \nabla),
$$
and hence (by functoriality) on
$$
\cL(\cE(\Einf(Y_n, \bu_0)))\ten \Cx \to \cL(\cE(\Gamma(\overline{Y}_{\!\!n},\widetilde{\bu_0}\ten \Omega_{\overline{Y}_{\!\!n}}^{*}(\log D_n)\ten \sA_{\Cx}^{0,*}))).
$$

Explicitly, we obtain these filtrations by considering the cofiltrations induced on the rings pro-representing these SDCs and DGLAs ($\cQ_*$ and $\cT_{\bullet}$ in the notation of \cite{sdc} Section \ref{sdc-sdcdgla}). The only technicality to check is that the resulting filtrations on the Lie algebra do indeed respect the Lie bracket (i.e. $[F^i,F^j] \subset F^{i+j}$). This is done by observing that we can recover the Lie bracket from the product (using a formula in some sense inverse to the Campbell-Baker-Hausdorff formula), giving a comultiplicativity condition on the cofiltrations satisfied by $\rho=*^{\sharp}$, and preserved by the functor $\cE$.

It then follows from \cite{Hodge3} that this is a mixed Hodge diagram in the sense of \cite{Morgan}, so we proceed as before.

\begin{theorem}\label{arbcxthm}
Let $S$ be the hull of the functor $\ddef_{\rho_0}$. Then 
$$
S \cong \R[[\H^1(X,\bu_0)^{\vee}]]/(f(\H^2(X,\bu_0)^{\vee})),
$$
where 
$$
f:\H^2(X,\bu_0)^{\vee} \to (\H^1(X,\bu_0)^{\vee})^2 \triangleleft \R[[\H^1(X,\bu_0)^{\vee}]]
$$ 
preserves the mixed Hodge structures.
In fact, the quotient map
$$
f:\H^2(X,\bu_0)^{\vee} \to (\H^1(X,\bu_0)^{\vee})^2/(\H^1(X,\bu_0)^{\vee})^3 \cong S^2(\H^1(X,\bu_0)^{\vee})
$$
is dual to half the cup product
\begin{eqnarray*}
\H^1(X,\bu_0) &\to&  \H^2(X,\bu_0),\\
a &\mapsto& \half a\cup a.
\end{eqnarray*}
\begin{proof} As shown above, $\ddef_{\rho_0}$ is governed by the SDC $\cE(\Einf(Y_n, \bu_0))^n$, hence by the DGLA $\cL(\cE(\Einf(Y_n, \bu_0)))$. The proof is now identical to Theorem \ref{morthm}. Note that since weight $0$ is permitted in both $\H^1(X,\bu)$ and $\H^2(X,\bu)$, we can draw no conclusions concerning the defining equations.
\end{proof}
\end{theorem}

\begin{remark}\label{subfields} 

As in \cite{Morgan}, we may replace $\R$ by any subfield $k$, taking $G/k$ algebraic. 

In fact, using Hodge theory for unitary local systems (as in \cite{Timm}), we may replace $\R$ by $\Cx$. As before, we obtain a bigraded decomposition on the complex cohomology, the only difference being that we must define the filtrations $F,\bar{F}$ separately, as complex conjugation no longer makes sense. In the general case, a sufficient hypothesis for $\bu_0$ to be a unitary local system becomes that 
$$
\ad(\rho_0)(\Gamma) \subset K,
$$
for some $K \le \GL(\g)$ a compact Lie group. 

Indeed, we may take any subfield  $k\subset \Cx$ instead of $\Cx$, giving a bigraded decomposition on complex cohomology, for which the corresponding weight decomposition will descend to $k$.
\end{remark}

\begin{remark} A much cleaner proof would be possible if the approach suggested in Remark \ref{hodgedef} were developed, since the filtrations on the SDCs
$$
 \cE(N_n)^n \to \cE(\Gamma(Y_n,\bu_0 \ten \sA^{\bullet}_{Y_n})\ten \Cx)^n \leftarrow \cE(\bar{N_n})^n, 
$$
where
$$
N_n=(\Gamma(\bar{Y_n},\widetilde{\bu_0}\ten \Omega_{\bar{Y_n}}^{*}(\log D_n)\ten \sA_{\Cx}^{0,*}), \nabla),
$$
allow us to define filtered deformation functors, with the same definitions as for DGLAs. This method would avoid having to make use of the noisome functor $\cL$.
\end{remark}

\section{Representations of the Algebraic Fundamental Group}\label{algrep}

Let $k=\bF_q$, fix a connected variety $X_k/k$, and let $X= X_k \ten_k \bar{k}$. Throughout, we will assume that $l$ is a prime not dividing $q$.

Fix a closed point $x$ of $X$, and denote the associated geometric point $x \to X$ by $\bar{x}$. 
The Weil group $W(X_k, \bar{x})$ is defined by:
$$
\xymatrix{
1 \ar[r] &\pi_1(X,\bar{x}) \ar[r] & \pi_1(X_k,\bar{x}) \ar[r] &\hat{\Z} \ar[r] & 0 \\
1 \ar[r]  & \pi_1(X,\bar{x}) \ar[r] \ar@{=}[u] &  W(X_k, \bar{x}) \ar@{^{(}->}[u]\ar[r] &\Z \ar@<-2pt>@{^{(}->}[u] \ar[r] & 0,
}
$$
with both rows exact. Each closed point $y \in |X_k|$  gives rise to a Frobenius element $\phi_y \in W(X_k,\bar{x})$ (as in \cite{Weil2} 1.1.8), defined up to conjugation, and  a representation of $W(X_k, \bar{x})$ is said to be pure of weight $n$ if the eigenvalues of all the Frobenius elements are algebraic, with all their complex conjugates of norm $q^{n/2}$.

 Both $\pi_1(X,\bar{x})$ and $W(X_k, \bar{x})$ are topological groups, the former with the profinite topology, and the latter with the product topology $W(X_k, \bar{x})\cong\pi_1(X_k,\bar{x})\by_{\hat{\Z}}\Z$.

\subsection{Representations to $\GL_n(\Q_l)$}

Fix a continuous representation 
$$
\rho_0:\pi_1(X,\bar{x}) \to \GL_n(\Q_l).
$$
We define the functor  
$$
\ddef_{\rho_0}:\C_{\Q_l} \to \Set
$$ 
by setting $\ddef_{\rho_0}(A)$ to be the set of isomorphism classes of continuous representations \mbox{$\rho:\pi_1(X,\bar{x}) \to \GL_n(A)$} deforming $\rho_0$.

We will first need to generalise the notion of constructible $\Q_l$-sheaf, as in \cite{Weil2} 1.1.1.

\begin{definition} For $A \in \C_{\Q_l}$, the category of constructible locally free $A$-sheaves will have as objects locally free $A$-sheaves $\sF$ which are  constructible  as $\Q_l$-sheaves in the sense of \cite{Weil2} 1.1.1. 
\end{definition}

It is well known that there is an equivalence of categories between continuous $\Q_l$-representations of the fundamental group, and constructible locally free $\Q_l$-sheaves. Let $\vv_0$ be the sheaf corresponding to $\rho_0$, with $\ww_0$ the underlying $\Z_l$-sheaf.

\begin{lemma}\label{vv}
For $A \in \C_{\Q_l}$, there is a functorial equivalence of groupoids between $\mathfrak{R}(A)$, the groupoid of continuous representations $\rho:\pi_1(X,\bar{x}) \to \GL_n(A)$ deforming $\rho_0$, and $\mathfrak{V}(A)$, the groupoid of rank $n$ constructible locally free $A$-sheaves deforming $\vv_0$.

\begin{proof}
The representation $\rho: \pi_1(X,\bar{x}) \to \GL_n(A)$ is continuous, so has compact image. Since $A$ is a finite dimensional  vector space over $\Q_l$, \cite{Se} LG 4 Appendix 1 provides the existence of a $\Z_l$-lattice $W$ generating $A^n$ such that $\rho$ factorises through $\GL(W)$. As in \cite{sdc} Section \ref{sdc-gln}, this gives a constructible $\Z_l$-sheaf $\ww$ on $X$. The $A$-module structure of $W\ten_{\Z_l} \Q_l$ now provides $\vv:=\ww\ten_{\Z_l} \Q_l$ with the structure of a constructible $A$-sheaf.

To construct a quasi-inverse, start with a constructible $A$-sheaf $\vv=\ww\ten_{\Z_l} \Q_l$. Let $W:=\ww_{\bar{x}}$ and $V=W \ten_{\Z_l} \Q_l$. Then, as in \cite{sdc} Section \ref{sdc-gln}, we obtain a representation 
$$
\rho:\pi_1(X,\bar{x}) \to \GL(W) \into \GL(V).
$$ 
We now choose an isomorphism $V \cong A^n$ compatible with the canonical isomorphism $(\vv_0)_{\bar{x}}\cong \Q_l^n$.
\end{proof}
\end{lemma}

\begin{definition}
For a constructible locally free $\Z_l$-sheaf $\sF$, define
$$
\CC^n(X,\sF):=\lim_{\substack{\longleftarrow \\m}} \CC^n(X,\sF/l^m\sF),
$$
where
$$
\CC^n(X,\sF):= \Gamma(X', (u^{*}u_*)^n u^{*}\sF),
$$
 for sheaves $\sF$ on $X$, and
$$
\sC^n(\sF):=(u_* u^{*})^{n+1}\sF,
$$
so that $\CC^n(X,\sF)=\Gamma(X, \sC^n(\sF))$.

For a constructible locally free $\Q_l$-sheaf $\sF \ten \Q_l$ of finite rank, define
\begin{eqnarray*}
\sC^n(\sF\ten_{\Z_l} \Q_l)&:=&\sC^n(\sF)\ten_{\Z_l} \Q_l,\\
\CC^n(X,\sF\ten_{\Z_l} \Q_l)&:=&\CC^n(X,\sF)\ten_{\Z_l} \Q_l.
\end{eqnarray*}
Note that this construction is independent of the choice of $\sF$, since $\sF$ is of finite rank.
\end{definition}

\begin{theorem}\label{defglnql}
Deformations of $\rho_0$ are described by the SDC
$$
E^n(A)=\exp(\CC^n(X,\ENd(\vv_0))\ten_{\Q_l}\m_A),
$$
with product given by the Alexander-Whitney cup product.

\begin{proof}
Given $\omega \in \mc_G(A)$, define 
$$
\vv_{\omega}:= \{ v \in \sC^0(\vv_0 \ten_{\Q_l} A): \alpha_{\sC^0(\vv_0 \ten A)}(v)=\omega \cup v\},
$$
where, for any sheaf $\sF$ on $X$, $\alpha_{\sF}$ is the canonical map $\alpha_{\sF}: \sF \to \sC^0(\sF)$.

Since $\m_A$ is a finite dimensional $\Q_l$ vector space, with nilpotent product, we may find a multiplicatively closed lattice $I \subset \m_A$ with
$$
\omega \in \exp(\CC^n(X,\ENd(\ww_0))\ten_{\Z_l} I).
$$
As in \cite{sdc} Section \ref{sdc-etshf}, this gives a constructible $(\Z_l\oplus I)$-sheaf $\ww_{\omega}$, and \mbox{$\vv_{\omega}=\ww_{\omega}\ten_{\Z_l} \Q_l$,} so is a constructible $A$-sheaf. 

Conversely, given a constructible $A$-sheaf $\vv$ deforming $\vv_0$, let $\ww$ be an underlying $\Z_l$-sheaf. We choose an isomorphism $\sC^0(\vv) \cong \sC^0(\vv_0)\ten_{\Q_l} \m_A$. The transition functions of $\ww$ then provide a suitable $\omega \in \mc_G(A)$. 
\end{proof}
\end{theorem}

\begin{remark} Equivalently, we may note that $L^n:=\CC^n(X,\ENd(\vv_0))$ has the structure of a DGAA (cup product being associative), with differential $d=\sum(-1)^i\pd^i$, so \emph{a fortiori} a DGLA. Given $\omega \in \mcl(A)$, let
$$
\vv_{\omega}:=\ker(d + \omega \cup:\sC^0(\vv_0)\ten_{\Q_l} A \to \sC^1(\vv_0)\ten_{\Q_l} A).
$$
That this gives the required deformation functor follows by observing that we have isomorphisms on tangent and obstruction spaces.
\end{remark}

\begin{lemma}\label{glnqlcoho}
If $\H^i(X,\ENd(\ww_0))$ and $\H^{i-1}(X,\ENd(\ww_0))$ are finitely generated, then 
$$
\H^i(E)=\H^i(X,\ENd(\vv_0)).
$$
\begin{proof}
\begin{eqnarray*}
\H^i(E)&=& \H^i(\CC^{\bullet}(X,\ENd(\vv_0)))\\
&=& \H^i(\CC^{\bullet}(X,\ENd(\ww_0)))\ten_{\Z_l}\Q_l\\
&=& \H^i( \lim_{\substack{\longleftarrow \\n}} \CC^{\bullet}(X,\ENd(\ww_0)\ten \Z/l^n))\ten_{\Z_l}\Q_l\\
&=& \H^i( \lim_{\substack{\longleftarrow \\_n}} \CC^{\bullet}(X,\ENd(\ww_0))\ten \Z/l^n)\ten_{\Z_l}\Q_l.
\end{eqnarray*}

Now, the tower
$$
\cdots \to \CC^{\bullet}(X,\ENd(\ww_0))\ten \Z/l^{n+1} \to \CC^{\bullet}(X,\ENd(\ww_0))\ten \Z/l^n \to \cdots
$$
clearly satisfies the Mittag-Leffler condition, so, by \cite{W} Theorem 3.5.8, we have the exact sequence:
$$
\xymatrix@=3ex{
0 \ar[r] & \Lim^1 \H^{i-1}(\CC^{\bullet}(X,\ENd(\ww_0))_n) \ar[r] \ar@{=}[d] & \H^i(\Lim_n \CC^{\bullet}(X,\ENd(\ww_0))_n) \ar[r] \ar@{=}[d]&  \Lim_n \H^i(\CC^{\bullet}(X,\ENd(\ww_0))_n) \ar[r] \ar@{=}[d]&    0\\
0 \ar[r]& \Lim^1 \H^{i-1}(X, \ENd(\ww_0)_n) \ar[r]& \H^i( \Lim_n \CC^{\bullet}(X,\ENd(\ww_0))_n) \ar[r]& \H^i(X,\ENd(\ww_0)) \ar[r]& 0.
}
$$
From \cite{Mi} Lemma V 1.11, it follows that the $\H^{i-1}(X, \ENd(\ww_0)\ten \Z/l^n)$ are finitely generated $\Z/l^n$-modules, so satisfy DCC on submodules, so this inverse system satisfies the Mittag-Leffler condition, making the $\Lim^1$ on the left vanish.

Thus
$$
\H^i(E)=\H^i(X,\ENd(\ww_0))\ten_{\Z_l}\Q_l= \H^i(X,\ENd(\vv_0)).
$$
\end{proof}
\end{lemma}

\subsection{Representations to an arbitrary group variety over $\Q_l$}\label{arbql}

Let $G$ be a  group variety over $\Q_l$, with associated Lie algebra $\g$. We will consider deformations over $\C_{\Q_l}$.
Given $A \in \C_{\Q_l}$, $G(A)$ has the structure of an $l$-adic Lie group.

Fix a continuous representation 
$$
\rho_0:\pi_1(X,\bar{x}) \to G(\Q_l).
$$
We define the functor  
$$
\ddef_{\rho_0}:\C_{\Q_l} \to \Set
$$ by setting $\ddef_{\rho_0}(A)$ to be the set of isomorphism classes of continuous representations \mbox{$\rho:\pi_1(X,\bar{x}) \to G(A)$} deforming $\rho_0$.

\begin{definition} Given a pro-$l$ group $K$, define a constructible principal $K$-sheaf to be a principal $K$-sheaf $\bD$, such that 
$$
\bD=\lim_{\substack{\longleftarrow \\{K \to F \text{ finite}}}} F \by^K \bD.
$$
Given an $l$-adic Lie group $G$, a constructible principal $G$-sheaf is a $G$-sheaf $\bB$ for which there exists a constructible principal $K$-sheaf $\bD$, for some $K \le G$ compact, with \mbox{$\bB=G \by^K \bD$} (observe that compact and totally disconnected is equivalent to pro-finite).
\end{definition}

\begin{lemma}
For $A \in \C_{\Q_l}$, there is a functorial equivalence of groupoids between $\mathfrak{R}(A)$, the groupoid of continuous representations $\rho:\pi_1(X,\bar{x}) \to G(A)$ deforming $\rho_0$, and $\mathfrak{B}(A)$, the groupoid of constructible principal $G(A)$-sheaves deforming $\Bu_0$.
\begin{proof}
Similar to Lemma \ref{vv}.
\end{proof}
\end{lemma}

Since $\ad \rho_0(\pi_1(X,\bar{x})) \le \GL(\g)$ is compact, the corresponding sheaf $\bu_0$ is a constructible $\Q_l$-sheaf of Lie algebras. Hence the sheaves $\sC^n(\bu_0)$ are sheaves of Lie algebras.

%
\begin{theorem}
Deformations of $\rho_0$ are described by the SDC
$$
E^n(A)=\exp(\CC^n(X,\bu_0)\ten_{\Q_l}\m_A),
$$
with product given by the Alexander-Whitney cup product. 
\begin{proof} As in Theorem \ref{defglnql} and \cite{sdc} Section \ref{sdc-smoothgroup}.
\end{proof}
\end{theorem}

\begin{lemma}
If we write $\bB=G\by^K\bD$, for some compact $K \le G$, and    $\H^i(X,\ad \bD_0)$ and $\H^{i-1}(X,\ad \bD_0)$ are finitely generated, then 
$$
\H^i(E)=\H^i(X,\ad \bB_0).
$$
\begin{proof} As for Lemma \ref{glnqlcoho}.
\end{proof}
\end{lemma}

\subsection{Structure of the hull}\label{weilhull}

 Given a mixed Weil sheaf $\sF$ over $\Q_l$ on $X$, for $X$ either smooth or proper, \cite{Weil2} shows that all the eigenvalues of Frobenius acting on the cohomology group $\H^i(X,\sF)$ are algebraic numbers $\alpha$, and for each $\alpha$, there exists a weight $n$, such that all complex conjugates of $\alpha$ have norm $q^{n/2}$. This provides us with a weight decomposition 
$$
\H^i(X,\sF)=\bigoplus_n \cW_n \H^i(X,\sF).
$$ 

\begin{theorem}\label{frobhull}
Let  $\rho_0: \pi_1(X,\bar{x}) \to G(\Q_l)$ be a continuous representation. Assume that
$$
\ad(\rho_0): \pi_1(X,\bar{x}) \to \Aut(\g),
$$
where $\Aut$ here denotes Lie algebra automorphisms, extends to a continuous mixed representation
$$
\widetilde{\ad(\rho_0)}: W(X_k,x) \to \Aut(\g).
$$
 Then the deformation functor
$$
\ddef_{\rho_0}:\C_{\Q_l} \to \Set
$$
has a hull of the form:
$$
\Q_l[[\H^1(X,\bu_0)^{\vee}]]/(f(\H^2(X,\bu_0)^{\vee})),
$$
where 
$$
f:\H^2(X,\bu_0)^{\vee} \to (\H^1(X,\bu_0)^{\vee})^2 \triangleleft \Q_l[[\H^1(X,\bu_0)^{\vee}]]
$$ 
preserves the Frobenius decompositions.
In fact, the quotient map
$$
f:\H^2(X,\bu_0)^{\vee} \to (\H^1(X,\bu_0)^{\vee})^2/(\H^1(X,\bu_0)^{\vee})^3 \cong S^2(\H^1(X,\bu_0)^{\vee})
$$
is dual to half the cup product
\begin{eqnarray*}
\H^1(X,\bu_0) &\to&  \H^2(X,\bu_0),\\
a &\mapsto& \half a\cup a.
\end{eqnarray*}
\begin{proof}
As in  \cite{Weil2}, the representation of the Weil group corresponds to an isomorphism $\Phi:F^*\bu_0 \to \bu_0$, where $F$ is the Frobenius endomorphism of $X$ over $\bF_q$. Since $\ad \rho_0$ is a representation to Lie algebra automorphisms, $\Phi$ is an isomorphism of constructible sheaves of Lie algebras. Hence we have the following morphisms of SDCs:
$$
\exp(\CC^n(X,\bu_0)\ten_{\Q_l}\m_A) \xra{F^*} \exp(\CC^n(X,F^*\bu_0)\ten_{\Q_l}\m_A) \xra{\CC^n(X,\Phi)} \exp(\CC^n(X,\bu_0)\ten_{\Q_l}\m_A).
$$

We therefore have a Frobenius endomorphism on our SDC. We could now proceed as in \cite{Weil2} \S 5, by forming the associated DGLA (as in \cite{sdc} Section \ref{sdc-sdctodgla}), which will, by functoriality, also have a Frobenius endomorphism. We would then form the minimal model, as in  Section \ref{smoothcx}, and construct a weight decomposition and quasi-isomorphism to a formal DGLA, as in \cite{Weil2} Corollary 5.3.7.

Instead, we will take a more direct approach, involving neither DGLAs nor minimal models. By abuse of notation, write $F^*$ for the Frobenius endomorphism on $E^{\bullet}$ given by the composition above. Since $F^*$ respects the structures, this gives us an endomorphism $F^*: \ddef_E \to \ddef_E$. Let $R$ be the hull of $\ddef_E$. We have:
$$
\xymatrix{
h_R \ar[r]^-{\etale} \ar@{-->}[d]_{\widetilde{F^*}} & \ddef_E  \ar[d]^{F^*}  \\
  h_R \ar[r]^-{\etale}& \ddef_E,  
}
$$
for some (non-unique) lift $\widetilde{F^*}$. Now let $\tilde{F}$ be the dual morphism $\tilde{F}:R \to R$ of complete local $\Q_l$-algebras. We have  canonical isomorphisms
$$
(\m_R/\m_R^2)^{\vee} \cong t_{\ddef_E} \cong \H^1(X,\bu_0),
$$
compatible with Frobenius. 

Since $\m_R/\m_R^n$ is a finite dimensional vector space over $\Q_l$ for all $n$, we may use the Jordan decomposition  (over $\bar{\Q_l}$) successively to lift the cotangent space $\m_R/\m_R^2$ to $\m_R/\m_R^n$, obtaining a space of generators $V \subset \m_R$, with the  map
$V \to \m_R/\m_R^2$ an isomorphism preserving the Frobenius decompositions.

Now look at the map 
$$
S=\Q_l[[V]] \onto R;
$$
let its kernel be $J \subset (V)^2$. $(J/VJ)^{\vee}$ is then a universal obstruction space for $h_R$. Since \mbox{$h_R \to \ddef_E$} is smooth, and $\H^2(X, \bu_0)$ is an obstruction space for $\ddef_{\rho_0}$, there is a unique injective map of obstruction spaces (\cite{Man} Proposition 2.18)
$$
(J/VJ)^{\vee} \into \H^2(X, \bu_0).
$$

Observe that both of these obstruction theories are  Frobenius equivariant: given  a functor $H$ on which $F$ (Frobenius) acts, we  say an obstruction theory $(W,w_e)$ is Frobenius equivariant if there is a Frobenius action on the space $W$ such that for every small extension
$$
e:\quad 0 \to I \to B \to A \to 0,
$$
and every $h \in H(A)$, we have $w_e(Fh)=Fw_e(h)$.

Since $h_R \to \ddef_E$ commutes with Frobenius, the corresponding map of obstruction spaces does (by uniqueness).

As before, we may use the Jordan decomposition on the successive quotients to obtain a space of generators $W \subset J$, such that $W \to J/VJ$ is an isomorphism preserving the  Frobenius decompositions.

Hence we have  maps
$$
\H^2(X, \bu_0)^{\vee} \onto  J/VJ \cong W \into \Q_l[[V]],
$$
and
$$
 V   \cong \m_R/\m_R^2  \cong \H^1(X, \bu_0)^{\vee},
$$
preserving the Frobenius decompositions.
Let $f:\H^2(X, \bu_0)^{\vee} \to \Q_l[[\H^1(X, \bu_0)^{\vee}]]$ be the composition.
Then
$$
R=\Q_l[[V]]/J \cong \Q_l[[\H^1(X, \bu_0)^{\vee}]]/(f(\H^2(X, \bu_0)^{\vee})),
$$
as required.

It only remains to prove the statement about the cup product, which follows since half the cup product is the primary obstruction map.
\end{proof}
\end{theorem}

\begin{corollary}
Let $X$ be smooth and proper, $\rho_0$ as above, and assume that
$
\widetilde{\ad(\rho_0)}
$
is pure of weight $0$. Then the deformation functor
$$
\ddef_{\rho_0}:\C_{\Q_l} \to \Set
$$
has a hull $R$ defined by homogeneous quadratic equations (given by the cup product). 
\begin{proof}
This follows since, under these hypotheses,  \cite{Weil2} Corollaries 3.3.4--3.3.6 imply that $\H^1(X,\bu_0)$ is pure of weight $1$, and $\H^2(X,\bu_0)$ is pure of weight $2$. 
\end{proof}
\end{corollary}

\begin{corollary}
Let $X$ be smooth, $\rho_0$ as above, and assume that
$
\widetilde{\ad(\rho_0)}
$
is pure of weight $0$. Then the deformation functor
$$
\ddef_{\rho_0}:\C_{\Q_l} \to \Set
$$
has a hull defined by equations of degree at most four.
\begin{proof}
This follows since, under these hypotheses,  \cite{Weil2} Corollaries 3.3.4--3.3.6 imply that $\H^1(X,\bu_0)$ is of weights $1$ and $2$, while $\H^2(X,\bu_0)$ is  of weights $2, 3$ and $4$. Hence the image of $f:\H^2(X,\bu_0)^{\vee} \to \Q_l[[\H^1(X,\bu_0)^{\vee}$ gives equations:
\begin{eqnarray*}
(\cW_{-1})^2 &\text{weight }-2 \\
(\cW_{-1})^3 +(\cW_{-1})(\cW_{-2}) &\text{weight }-3\\
(\cW_{-1})^4 +(\cW_{-1})^2(\cW_{-2}) +(\cW_{-2})^2 &\text{weight }-4,
\end{eqnarray*}
hence equations of degree at most $4$.
\end{proof}
\end{corollary}

\begin{remarks}\label{qlbar}
\begin{enumerate}
\item
In the above working, we may replace $\Q_l$ by any finite field extension $E/\Q_l$. Indeed, since $\pi_1(X,\bar{x})$ is pro-finitely generated, any representation to $G(\bar{\Q_l})$ will factor through $G(E)$ for some finite extension $E$ of $\Q_l$, so we may replace $\Q_l$ by $\bar{\Q_l}$, with the usual conventions for constructible $\bar{\Q_l}$-sheaves. 

\item
If $G=\GL_n$, a large class of examples can be produced to satisfy the hypotheses. For each irreducible continuous representation $\tilde{\rho_0}: W(X_k,\bar{x}) \to \GL_n(E)$, the induced representation $\rho_0: \pi_1(X,\bar{x}) \to \GL_n(E)$ satisfies the hypotheses. This follows from \cite{Weil2} Conjecture 1.2.10 (and remarks preceding) which imply that every irreducible Weil sheaf over $\bar{\Q_l}$ is of the form
$$
\sF = \sP \ten \bar{\Q_l}^{(b)},
$$
for a  pure sheaf $\sP$. Here, $\bar{\Q_l}^{(b)}$ is the constant sheaf $\Q_l$, on which the Frobenius action is multiplication by $b$,  for $b \in \bar{\Q_l}^*$.

Hence
$$
\ENd(\sF)\cong \sF \ten \sF^{\vee} \cong \sP \ten \sP^{\vee},
$$
which is pure of weight $0$. Setting $\widetilde{\ad \rho_0}:= \ad \tilde{\rho_0}$ gives the desired result.

Lafforgue has proved \cite{Weil2} Conjecture 1.2.10 (known as the Deligne Conjecture) as \cite{La} Theorem VII.6.
\end{enumerate}
\end{remarks}
\bibliographystyle{plain}
\addcontentsline{toc}{section}{Bibliography}
\bibliography{references.bib}
\end{document}